\documentclass[12pt,twoside]{article}
\usepackage[all]{xy}
\usepackage{amssymb}
\begin{document}

\font\BBb=msbm10 at 12pt
\newcommand{\CC}{{\hbox{\BBb{C}}}}
\newcommand{\QQ}{{\hbox{\BBb{Q}}}}
\newcommand{\FF}{{\hbox{\BBb{F}}}}
\newcommand{\NN}{{\hbox{\BBb{N}}}}
\newcommand{\RR}{{\hbox{\BBb{R}}}}
\newcommand{\ZZ}{{\hbox{\BBb{Z}}}}
\newcommand{\MM}{{\hbox{\BBb{M}}}}
\newcommand{\BB}{{\hbox{\BBb{B}}}}
\newcommand{\PP}{{\hbox{\BBb{P}}}}
\newcommand{\DD}{{\hbox{\BBb{D}}}}
\newcommand{\TT}{{\hbox{\BBb{T}}}}

\newcommand{\arr}{\longrightarrow}
\newcommand{\cK}{{\widehat K}^{(q)}}
\newcommand{\s}{\subset}
\newcommand{\sm}{\setminus}
\newcommand{\f}{\varphi}
\newcommand{\e}{\epsilon}
\newcommand{\fa}{\textfrak{a}}
\newcommand{\fb}{\textfrak{b}}

\newcommand{\reco}{\mbox{$\subset\hspace{-.133 in}\raisebox{.25ex}
{$\scriptscriptstyle\subset$}$}}

\newcommand{\fin}{\hspace*{\fill}\rule{2.8mm}{2.8mm}}

\newtheorem{thm}{Theorem}
\newtheorem{lem}{Lemma}
\newtheorem{ex}{Example}
\newtheorem{re}{Remark}
\newtheorem{co}{Corollary}
\newtheorem{pp}{Proposition}
\newtheorem{df}{Definition}
 
\pagestyle{myheadings}

\markboth{V\^{\i}j\^{\i}itu}{Rad\'{o} 
theorem for complex spaces} 

\title{A Rad\'{o} theorem for complex spaces} 
\author{Viorel V\^{\i}j\^{\i}itu }
\date{}
\maketitle
\noindent {\it Abstract}. We generalize Rad\'{o}'s extension theorem 
from the complex plane to reduced complex spaces.
\vspace{0.2cm}

\vspace{0.2cm}


\section{Introduction} 
A theorem due to Rad\'{o} asserts that
a continuous complex valued function 
on an open subset of the complex plane 
is holomorphic 
provided that it is holomorphic off its zero set.

Essentially this theorem was proved in \cite{Ra-24}. 
Since then many other proofs have been proposed, {\it e.g.} \cite{BeSt-51}, \cite{Ca-52}, \cite{He-56}, and \cite{Ka-67}. 
The articles \cite{Au-51}, \cite{Ri-83} and \cite{St-68}
give some generalizations.

Rad\'{o}'s statement remains true for complex manifolds 
(or, more generally, for normal complex spaces)
as well as in the complex plane. 

In this short note we investigate a natural 
extension of Rad\'{o}'s theorem when the ambient space
has (non normal) singularities. 

Complex spaces, unless explicitely stated,
are assumed to be reduced and countable at infinity.
Let $\NN=\{1,2, \ldots \}$.

Here we state our main results. 
\begin{pp} There is an irreducible Stein curve 
$X$ and a continuous function $f:X \arr \CC$
that is holomorphic off its zero set but no 
power $f^n$, $n \in \NN$, is globally holomorphic.
\end{pp}
\begin{thm} Let $X$ be a complex space
and $\Omega \s X$ be a relatively compact open set. 
Then, there is $n_\Omega \in \NN $ such that, 
for every continuous function $f:X \longrightarrow \CC$ 
that is holomorphic off its zero set, 
for every integer $n \geq n_\Omega$, the power
$f^n$ is holomorphic on $\Omega$.
\end{thm}

Recall the following definition. Let $X$ be a complex space.
A function $f:X \arr \CC$ is {\it c}-{\it holomorphic}
if it is continuous and the restriction of $f$ to the subset
${\rm Reg}(X)$ of manifold points of $X$ 
is holomorphic \cite{Wh-72}.
\section{Proof of Proposition 1}
The singularities that we implant at the points
$2,3, \ldots$, of $\CC$ are generalized cusp singularities. 

Let $p$ and $q$ be integers $\geq 2$ that are coprime. 
Consider the cusp like irreducible and locally irreducible complex curve
$$\Gamma=\{(z_1,z_2) \in \CC^2 \;:\, z_1^p=z_2^q\} 
\s \CC^2.$$ 
Its normalization is $\CC$ and 
$\pi:\CC \longrightarrow \Gamma $, $t \mapsto (t^q,t^p)$,
is the normalization map. Note that $\pi$ is a
homeomorphism.

A continuous function 
$h:\Gamma \arr \CC$ that is holomorphic off its zero set, 
but fails to be globally holomorphic
is produced as follows. 

Select natural numbers $m$ and $n$ 
with $mq-np=1$, and define 
$h:\Gamma \longrightarrow \CC$ by setting 
for $(z_1,z_2) \in \Gamma$,
\[ h(z_1,z_2):= \left\{
\begin{array}{ll}
{z_1^m}/{z_2^n} & \mbox{ if } z_2 \neq 0, \\   
0 & \mbox{ if } z_2=0. \\
\end{array}
\right.
\]
It is easily seen that $h$ is continuous (as $\pi$ is a homeomorphism, the continuity of $h$ follows from that
of $h \circ \pi$, that equals ${\rm id}$ of $\CC$), $h$ is holomorphic off its zero set 
(incidentally, here, this is the set of smooth points of
$\Gamma$), and $h$ is not holomorphic about $(0,0)$ (use a Taylor series expansion
about $(0,0) \in \CC^2$ of a presumably holomorphic extension). 

Furthermore, $h^k$ is globally holomorphic provided that 
$k \geq (p-1)(q-1)$. (Because every integer
$ \geq (p-1)(q-1)$ can be written 
in the form $\alpha p + \beta q$ with 
$\alpha,\beta \in \{0, 1, 2, \ldots \}$,
and since $h^p$ and  $h^q$ are holomorphic being the restrictions of $z_2$ and $z_1$ to $\Gamma$ respectively.)

Also, $z_1^a z_2^b h$ is holomorphic on $\Gamma$ provided that
$$ q \lfloor(m+a)/p \rfloor+b \geq n,$$
where $\lfloor \cdot \rfloor$ is the floor function.

It is perhaps interesting to note that 
the stalk of germs of weakly holomorphic
functions $\widetilde{\cal O}_0$ at $0$
is generated as an ${\cal O}_0$-module by the 
germs at $0$ of $1, h, \ldots, h^{r}$, where 
$r=-1+\min\{p,q\}$..
\vspace{0.2cm}

Now, for each integer $k \geq 2$, let 
$\Gamma_k:=\{(z_1,z_2) \in \CC^2 \;; \, 
z_1^k=z_2^{k+1}\}$. As previously noted, 
$\Gamma_k$ is an irreducible curve 
whose normalization map is $\pi_k:\CC \longrightarrow \Gamma_k, \, t \mapsto (t^{k+1},t^k)$, and the function 
$h_k:\Gamma_k \longrightarrow \CC$ defined for $(z_1,z_2) 
\in \Gamma_k$ by
\[ h(z_1,z_2):= \left\{
\begin{array}{ll}
{z_1}/{z_2} & \mbox{ if } z_2 \neq 0, \\   
0 & \mbox{ if } z_2=0. \\
\end{array}
\right.
\]
has the following properties:
\begin{enumerate}
\item[${\bf a}_k)$]{The function $h_k$ is weakly holomorphic.}
\vspace{-0.2cm}

\item[${\bf b}_k)$]{The power $h_k^{k-1}$ is not holomorphic.}

\item[${\bf c}_k)$]{The function $z_1^{k-1}h_k$ is holomorphic
because it is the restriction of $z_2^k$ to $\Gamma_k$.}
\end{enumerate}

Here, with these examples of singularities at hand, 
we change the standard complex structure of 
$\CC$ at the discrete analytic set 
$\{2,3, \ldots \}$ by complex surgery in order to obtain an irreducible Stein complex curve $X$ and a discrete subset 
$\Lambda=\{x_k: k=2, 3, \ldots \}$ such that, 
at the level of germs
$(X,x_k)$ is biholomorphic to $(\Gamma_k,0)$. 

The surgery, that we recall 
for the commodity of the reader because in some 
monographs like \cite{KK-83}
the subsequent condition $(\star)$ is missing, goes as follows,

Let $Y$ and $U'$ be complex spaces together with 
analytic subsets $A$ and $A'$ of $Y$ and $U'$ respectively, 
such that there is an open neighborhood $U$ of $A$ in $Y$ and
$\f:U \sm A \arr U' \sm A'$ that is biholomophic.

Then define 
$$X:=(Y \sm A) \sqcup_\f U':= (Y \sm A) \sqcup U'/_{\sim}$$
by means of the equivalence relation 
$U \sm A \ni y \sim \f(y) \in U' \sm A'$. 

Then there exists exactly one complex structure 
on $X$ such that $U'$ and 
$Y \sm A$ can be viewed as open subsets of $X$ 
in a canonical way provided that the 
following condition is satisfied:
\begin{enumerate}
\item[$(\star)$] {For every $y \in \partial U$ and $a' \in A'$ there are open neighborhoods 
$D$ of $y$ in $Y$, $D \cap A=\emptyset$, and 
$B$ of $a'$ in $U'$ such that
$\f(D \cap U) \cap B \subseteq A'.$}
\end{enumerate}
Thus $X$ is formed from $Y$ by ''replacing`` $A$ with $A'$.

In practice, the condition $(\star)$ is fulfilled if $\f^{-1}: U' \sm A' \arr U \sm A$ 
extends to a continuous function $\psi:U' \arr U$ such that $\psi(A')=A$. In this case, if $D$ and
$V$ are disjoint open neighborhoods of $\partial U$ and $A$ in $Y$ respectively, then 
$B=A' \cup \f(V \sm A)$ is open in $U'$ because it equals $\psi^{-1}(V)$ and  $(\star)$ follows immediately.
(This process is employed, for instance,
in the construction of the blow-up of a point in a complex manifold!)

Coming back to our setting, consider $Y=\CC$, $A=\{2,3, \ldots \}$ and for each $k=2,3, \ldots$,
let $\Delta(k,1/3)$ be the disk in $\CC$ centered at $k $ of radius $1/3$ that is mapped holomorphically
onto an open neighborhood $U_k$ of $(0,0) \in \Gamma_k$ through the holomorphic map $t \mapsto \pi_k(z-k)$.
Applying surgery, we get an irreducible Stein curve $X$ and the discrete subset $\Lambda $ with the aforementioned 
properties. 

It remains to produce the function $f$ as stated. For this
we let ${\cal I} \s {\cal O}_X$ be the 
coherent ideal sheaf with support $\Lambda$
and such that ${\cal I}_{x_k}= \mathfrak{m}^{k-1}_{x_k}$
for $k=2,3 \ldots$, where $\mathfrak{m}_{x_k}$
is the maximal ideal of the analytic algebra of the stalk of ${\cal O}_X$ at $x_k$.

From the exact sequence
$$ 0 \arr {\cal I} \arr \widetilde{\cal O} \arr \widetilde{\cal O}/{{\cal I}} \arr 0, $$
where $\widetilde{\cal O}$ stands for the sheaf of germs of weakly holomorphic functions
in $X$, 
we obtain a weakly holomorphic function $f$ on $X$ 
such that, for each $k=2,3, \ldots$,
at germs level $f$ equals $f_k$ (mod ${\cal I}_{x_k}$). 

By  properties ${\bf a}_k)$, ${\bf b}_k)$ and 
${\bf a}_k)$ from above, it follows that
no $n \in \NN$ exists for which $f^n$ becomes 
holomorphic on $X$. (For instance, if $f=f_k+g_k^{k-1}$, 
for certain $g_k \in \mathfrak{m}_{x_k}$, then
$f^{k-1}$ is not holomorphic about $x_k$.)
\section{Proof of Theorem 1}
Before starting the proof, we collect a few definitions and
facts.

Let $(A,a)$ be a germ of a local analytic subset in $\CC^n$.
The third Whitney cone ${\rm C}_a(A)$  
consists of all vectors ${\bf v} \in \CC^n$ such that there
is a sequence $(a_j)_j$ of points in $A$ converging to $a$ and a sequence $(n_j)_j$ of $\NN$ such that
$$ n_j(a_j-a) \mapsto {\bf v} \; \mbox{   as  } \:
j \mapsto \infty.$$
If $A$ is pure $k$-dimensional, then ${\rm C}_a(A)$ is also a pure $k$-dimensional
analytic set. This cone ${\rm C}_a(A)$ lies in the tangent space $T_a(A)$ to $A$ at $a$ 
($T_a(A)$ is the sixth Withney cone; see \cite{Ch-89}
for more details).

Assume that the analytic germ $(A,0) \s \CC^k_z \times \CC_w^{m-k}$,
with $x=(z,w)$, is pure $k$-dimensional such that the projection 
$\pi(z,w)=z$ is proper on it, and $\pi^{-1}(0) \cap {\rm C}_0(A)=\{0\}$. 
Then $\pi$ is a branched covering on $A$ with covering number $d:=\mbox{deg}_0 A$
and critical set $\Sigma$. Now, whenever $x \not \in \Sigma$, we may define, for a non-constant ${\rm c}$-holomorphic
germ $f:(A,0) \longrightarrow (\CC,0)$, the polynomial
$$\omega(x,t)= \prod_{\pi(x')=\pi(x)} (t-f(x')).$$
Since $f$ is holomorphic on the regular part $\mbox{Reg}(A)$ of $A$ and $f$ is continuous on $A$, 
we obtain a distinguished Weierstrass polynomial of degree $d$ with holomorphic coefficients,
$W(z,t)= t^d + a_1(z) t^{d-1}+ \cdots + a_d(z)$, such that $W(x,f(x))=0$ for $x \in A$.

Note that $W $ does not depend on $w$. Besides, if $W(t,z)=0$, then
$$|t|={\rm O}(\|z\|^{1/d}) \mbox{  as } (z,t) \rightarrow 0$$
meaning that there are positive constants $M$ and $\e $ such that, if $W(z,t)=0$ and 
$\max\{|t|, \|z\|\}<\e$, then  $|t| \leq M \|z\|^{1/d}$. 
\vspace{0.2cm}

Therefore, if  $f:(A,0) \longrightarrow (\CC,0)$ is a non-constant c-holomorphic germ on a pure $k$-dimensional
analytic germ at $0 \in \CC^m$, then 
$$|f(\zeta)|=O(\|\zeta\|^{1/d})  \mbox{ as  } \zeta \rightarrow 0,$$ 
where $d$ equals the degree of the covering of the projection $\pi$ satisfying $\pi^{-1}(0) \cap {\rm C}_0(A)=\{0\}$.
\vspace{0.2cm}

In other words, if $f:(A,a) \longrightarrow (\CC,0)$ is a non-constant c-holomorphic germ on a pure dimensional analytic germ $A$ at $a$, we may define the {\it order of flatness}
of $f$ at $a$ to be 
$$\mbox{\rm ord}_a f= \max \{ \alpha>0 \,:\,  |f(x)| ={\rm O}(\|x-a\|^\alpha) \mbox{ as } x \rightarrow a \}.$$
It follows from the above discussion  that $\mbox{\rm ord}_a f \geq 1/d$.
\vspace{0.2cm}

Following Spallek \cite{Sp-65}, a germ function $f:(A,a) \longrightarrow (\CC,0)$ 
is said to be ${\rm O}^N$-{\it approximable at} 
$a$ if there exists a polynomial $P(z,\bar z) $ of degree at most $N-1$ 
in the variables $z_i-a_i, \overline{z_i-a_i}$, $i=1, \ldots, n$, such that
$$ |f(z)-P(z,\bar z)|={\rm O}(\|z-a\|^N) \mbox{ as } z \rightarrow a. $$

Here we quote from Siu \cite{Siu-69} the following
result that improves onto Spallek's similar one 
({\it loc. cit.}).
\begin{pp} For every compact set $K$ of a complex space $X$ there exists a positive integer 
$N$ $($depending on $K$$)$ such that, if $f$ is a 
{\rm c}-holomorphic function germ
at $x \in K$ and the real part $ \Re f$ of $f$ is ${\rm O}^N$-approximable in some neighborhood of $x$, then $f$ is 
a holomorphic germ at $x$.
\end{pp}
The above discussion concludes readily the proof of Theorem 1. 
\section{A final remark}
Below we answer a question raised by Th. Peternell at the XXIV Conference on
Complex Analysis and Geometry, held in Levico-Terme, June 10--14, 2019. He asked whether or not a similar statement 
like Theorem 1 does hold for non reduced complex spaces.

More specifically, let $(X,{\cal O}_X)$ be a not necessarily reduced complex space
and 
$f:X \arr \CC$ be continuous such that, if $A$ denotes the zero set of $f$, then $X \sm A$ is dense in $X$ and
there is a section $\sigma \in \Gamma(X \sm A,{\cal O}_X)$ whose reduction 
${\rm Red}(\sigma)$ equals $f|_{X \sm A}$. 

Is it true that, for every relatively compact open subset
$D$ of $X$, there is a positive integer $n$ such that $\sigma^n$ extends to a section
in $\Gamma(D,{\cal O}_X)$ ?
\vspace{0.2cm}

We show that the answer is ''No``. 

Recall that, if $R$ is a commutative ring with unit 
and $M$ is an $R$-module, we can endow
the direct sum $R \oplus M$ with a ring structure with the obvious addition, and multiplication
defined by 
$$ (r,m) \cdot (r',m')=(rr',rm'+r'm).  $$
This is the Nagata ring structure from algebra \cite{Na-62}. 

Now, if $(X,{\cal O}_X)$ is a complex space, and $\cal F$ a coherent ${\cal O}_X$-module, 
then ${\cal H}:={\cal O}_X \oplus {\cal F}$ 
becomes a coherent sheaf of analytic algebras and
$(X,{\cal H})$ a complex space (\cite{Fo-67}, Satz 2.3). 

The example is as follows. Let ${}_n{\cal O}$ denotes the structural sheaf of $\CC^n$. The above discussion
produces a complex space  $(\CC,{\cal H})$ such that 
${\cal H}={}_1{\cal O} \oplus {}_1{\cal O}$,
that can be written in a suggestive way 
${\cal H}={}_1{\cal O}+\e \cdot {}_1{\cal O}$, 
where $\e$ is a symbol with $\e^2=0$. 
As a matter of fact, if we consider $\CC^2$ with complex 
coordinates $(z,w)$ and the
coherent ideal ${\cal I}$ generated by $w^2$, then ${\cal H}$
is the analytic restriction of the quotient 
${}_2{\cal O}/{{\cal I}}$ to $\CC$. 

The reduction of $(\CC,{\cal H})$ is $(\CC,{}_1{\cal O})$. 
A holomorphic section of ${\cal H}$ over an open set 
$U \s \CC$ consists of couple of ordinary
holomorphic functions on $U$. 

Now take $f$ the identity function ${\rm id}$ on $\CC$, and 
the holomorphic section 
$\sigma \in \Gamma(\CC^\star,{\cal H})$ 
given by $\sigma= {\rm id} + \e g$,
where $g$ is holomorphic on $\CC^\star$ 
having a singularity at $0$, for
instance $g(z)= 1/z$. 

Obviously, the reduction of $\sigma$ is 
the restriction of ${\rm id}$
on $\CC^\star$, and 
no power $\sigma^k$ of $\sigma$
extends across $0$ to a section in 
$\Gamma(\CC, {\cal H})$, since $\sigma^k={\rm id} + \e kg$ and
$g$ does not extend holomorphically across $0 \in \CC$.

\vspace{0.5cm}

\noindent  Universit\'{e} de Lille, Lab. Paul Painlev\'{e}, 
B\^{a}t. M2

\noindent F-59655 Villeneuve d'Ascq Cedex, France

\noindent E-mail: {\tt viorel.vajaitu@univ-lille.fr}
\end{document}